\documentclass[11pt]{article}
\usepackage[english]{babel}
\usepackage{amsfonts, amsmath, amsthm, amssymb,amscd,indentfirst}
\usepackage[notref,notcite]{showkeys}
\usepackage{graphicx}
\usepackage{wrapfig}
\usepackage[rflt]{floatflt}
\usepackage{amsmath,amssymb,latexsym,indentfirst}
\usepackage{enumerate,pst-all}
\usepackage{makeidx}
\usepackage{color}
\usepackage{times}
\usepackage{palatino}

\usepackage[latin1]{inputenc}
\usepackage{amstext}
\usepackage{amsfonts}
\newtheorem{thm}{Theorem}[section]
\newtheorem{defi}{Definition}[section]
\newtheorem{lema}{Lemma}[section]
\newtheorem{cor}{Corollary}[section]

\newcommand{\s}{\sigma}
\newcommand{\cc}{{\cal C}}
\newcommand{\om}{\Omega}

\newcommand{\br}{\Bbb R}
\newcommand{\bn}{\Bbb N}

\newcommand{\W}{{\cal W}}

\renewcommand{\div}{{\rm div \, }}
\author{Barnab\'e Pessoa Lima \and J. F\'abio  Montenegro \and  Newton Lu\'{\i}s Santos}
\title{Eigenvalues Estimates for the $p$-Laplace Operator on Manifolds}
\begin{document}

\author{\and Barnab\'e P. Lima\thanks{Universidade Federal do
Piau\'{\i},
 Departamento de Matem\'atica, Campus Petronio Portela, Ininga, 64049-550 Teresina/PI, Brazil ({\tt barnabe@ufpi.br}). Partially supported by Instituto do Mil\^enio} \and J. F\'abio  Montenegro \thanks{Universidade Federal do Cear\'a,
Departamento de Matem\'atica, Campus do Pici, R. Humberto Monte, s/n, 60455-760,
Fortaleza/CE, Brazil ({\tt fabio@mat.ufc.br}). Partially supported by CNPq, grant ...}
 \and Newton L. Santos\thanks{Universidade Federal do
Piau\'{\i},
 Departamento de Matem\'atica, Campus Petronio Portela, Ininga, 64049-550 Teresina/PI, Brazil ({\tt newtonls@ufpi.br}). Partially supported by a CNPq Posdoctoral Grant  155723/2006-5 }   } \maketitle

\begin{center} Dedicated to Prof. Manfredo P. do Carmo in his 80th birthday\end{center}

\begin{abstract}
\noindent We obtain geometric estimates for the first eigenvalue and the fundamental tone of the $p$-laplacian on manifolds in terms of admissible vector fields. Also, we defined a new spectral invariant and we show its relation with the geometry of the manifold.
\end{abstract}

\section{Introduction and statement of results}

The Laplace-Beltrami operator on a Riemannian manifold, its spectral theory and the relations between its first eigenvalue and the geometrical data of the manifold, such as curvatures, diameter, injectivity radius, volume, has been extensively studied in the recent mathematical literature. In the last few years, another operator, called $p$-Laplacian, arising from problems on  Non-Newtonian Fluids, Glaceology, Nonlinear Elasticity, and  in problems of Nonlinear Partial Differential Equations came to the light of Geometry. Since then, geometers showed that this  operator exhibit some very interesting analogies with the Laplacian.

Let $(M,g)$ be a smooth Riemannian manifold and $\om\subset M$ a domain. For $1< p< \infty$, the $p$-laplacian on $\om$ is defined by
\begin{equation}
\triangle_p(u)=-\div \Bigl[\|\nabla u\|^{p-2} (\nabla u)\Bigr].
\end{equation}
This operator appears naturally from the variational problem associated to the energy functional
\begin{equation*}
E_p: W^{1,p}_0(\om)\to \br\qquad\mbox{given by}\qquad E_p(u)=\int_\om \|\nabla u\|^p\,\,\, d\om
 \end{equation*}
  where $W^{1,p}_0(\om)$ denotes the Sobolev space given by the closure of ${\cal C}^{\infty}(\om)$-functions with  compact support in $\om$ for the norm
  $$\|u\|_{1,p}^p=\int_\om |u|^p\,\,\, d\om+\int_\om \|\nabla u\|^p\,\,\, d\om.$$
Observe that, when $p=2$, $\triangle_2$ is just the Laplace-Beltrami operator.
We are interested in the nonlinear eigenvalue problem
\begin{equation}
\triangle_p u+\lambda |u|^{p-2}u=0 \label{eigenv problem}
\end{equation}
Since solutions for this problem, for arbitrary  $p\in(1,\infty)$ are only locally ${\cal C}^{1,\alpha}(\om)$ (exceptions for the case $p=2$), they must be described in the sense of distributions, that is, $u\in W^{1,p}_0(\om)\setminus \{0\}$ is an eigenfunction associated to the eigenvalue $\lambda$, if
$$
\int_\om \|\nabla u\|^{p-2}g(\nabla u,\nabla \phi)\,\,\, d\om =\lambda \int_\om |u|^{p-2} u\phi \,\,\, d\om
$$
for any test function $\phi\in \cc^{\infty}_0(\om)$. Clearly, $\lambda=0$ is an eigenvalue of $\triangle_p$, with associated  eigenfunctions being the constant functions. 
 \\
 Question: Does there exists spectrum for $\triangle_{p}$  $\lambda=0$ is an eigenvalue, but are there others ? This should be clear in this text.

The set $\s_{p}(M)$ of the remaining eigenvalues of (\ref{eigenvproblem}) is an unbounded subset of $(0,\infty)$ at least for euclidian domains, $\om\subset \br^n$, as quoted in \cite{Matei1} (see also \cite{Lindqvist}) whose infimum $\inf \s_p=\mu_{1,p}(\om)$ is an eigenvalue.

 Question:  What about for manifolds?  Is the set $\s_{p}(M)$ unbounded? non-empty etc
 \\
 
 It is also known (see \cite{Matei1}) that the first eigenvalue is simple and the first eigenfunctions for geodesic balls on space-forms are radial.
\\

 From here on, things are not clear at all!!! Is there a Rayleigh Theorem for  $\triangle_{p}$  Is the fundamental tone of an open domain the first eigenvalue

Let $\om\subset M$ a domain with compact closure and nonempty boundary $\partial \om$. The {\it $p$-fundamental tone of $\om$}, denoted by $\mu_{p}^*(\om)$ is defined as follows:
$$\mu_{p}^*(\om)=\inf\left\{ \frac{\int_\om \|\nabla f\|^p\,\,\, d\om}{\int_\om |f|^p\,\,\, d\om}; f\in W_0^{1,p}(\om), \,\,\, f\neq 0\right\}$$
when $ \om$ has piecewise smooth boundary, then $\mu^*_p(\om)$ coincides with the first eigenvalue of the  problem (\ref{eigenv problem}) with boundary condition $u\vert_{\partial\om}= 0$, by Rayleigh's Theorem.
In particular, when $M$ is a closed manifold, i.e. compact without boundary, we get
$$\mu_p^*(M)=\mu_{1,p}(M)=\inf\left\{ \frac{\int_M \|\nabla u\|^p\,\,\, dM}{\int_M |u|^p\,\,\, dM}; u\in W^{1,p}(M), u\neq 0, \int_M|u|^{p-1}u\,\,\, dM=0 \right\}$$
Observe that if $\om_1\subset \om_2$ are bounded domains, then $\mu_{p}^{\*}(\om_1)\ge \mu_{p}^{\*}(\om_2)\ge 0$. Thus one may define the $p$-fundamental tone $\mu_{p}^{\*}(M)$ of an open Riemannian manifold (i.e., complete noncompact) as the limit
$$\mu_{p}^*(M)=\lim_{r\to\infty}\mu_{p}^*(B_r(q))$$
where $B_r(q)$ is the geodesic ball of radius $r$ centered at $q\in M$

When $p=2$ the $p$-laplacian is simply the laplacian and the $p$-fundamental tone is simply called the fundamental tone. Interesting estimates on the fundamental tone for the Laplace-Beltrami operator on a Riemannian manifold have been recently obtained by G. Pacelli Bessa and the second author (see, for instance \cite{Bessa-Montenegro1} and \cite{Bessa-Montenegro2}). This paper presents an attempt to extend their variational argument to the $p$-laplacian.
Precisely we have
\begin{thm}\label{main1}
Let $\om\subset M$ be a domain ($\partial \om\neq \emptyset$) in a Riemannian manifold, $M$. Then
\begin{equation}
\mu_{p}^*(\om)\ge \frac{c(\om)^p}{p^p}>0
\end{equation}
where $c(\om)$ is the constant given in (\ref{fundamental constant})
\end{thm}

To present the second variational estimate we need to introduce some preliminary definitions which will allow us to deal with divergence of vector fields in a weak sense.

\vspace{.5cm}

\begin{defi}[Weak divergence] Let $(M,g)$ be a Riemannian manifold and $X\in L^1_{loc}(\mathfrak{X}(M))$ (in the sense that $\|X\|\in L^1_{loc}(M)$) A function $h\in L^1_{loc}(M)$ is said to be a weak divergence of $X$, denoted by $h={\rm Div}X$ if for every $\phi\in \cc^{\infty}_0(M)$ it holds
\begin{equation}
\int_M \phi h\,\, d\mu=-\int_Mg(\nabla \phi, X)\,\, d\mu
\end{equation}
\end{defi}

\vspace{.5cm}

The weak divergence exists for almost every point of $M$. If $X\in {\cal W}^{1,1}(M)$ and $f\in\cc^1(M)$ then $fX\in {\cal W}^{1,1}(M)$ with ${\rm Div} (fX)=g(\nabla f,X)+f{\rm Div} X$. In particular for $f\in\cc^{\infty}_0(M)$ we have that
\begin{equation}
\int_M {\rm Div}(fX)\,\, d\mu=\int_M \Big[g(\nabla f,X)+f{\rm Div}(X)\Big]\,\, d\mu =0
\end{equation}
With these notations fixed we have

\vspace{.5cm}

\begin{thm}\label{main2}
Let $(M,g)$ be a Riemannian manifold. Then the following estimate holds
\begin{equation}
\mu_{p}^*(M)\ge \sup\left\{ \inf_\om\Big((1-p)\|X\|^q+{\rm Div}(X)\Big), \,\,\, X\in \W^{1,1}(M)\right\}
\end{equation}
\end{thm}

As a simple consequence of Theorem \ref{main1} we obtain a generalization of Mckean  theorem, \cite{Mckean}

\vspace{.5cm}

\begin{thm}[Generalized Mckean]\label{mckean}
Let $M$ be an $n$-dimensional, complete noncompact, simply connected Riemannian manifold with sectional curvature $K\leq -c^2<0$, then
$$
\mu_p^*(M)\ge \frac{(n-1)^pc^p}{p^p}
$$
In particular, when $p=2$ this is the Mckean theorem.
\end{thm}

\vspace{.5cm}

Contrary to the Laplace operator, the $p$-laplacian has not been proved to be discrete, even for euclidian domains $\om \subset \br^n$ (as remarked in \cite{Idrissa}). There are few results related to the spectrum of such operator. For instance, Lindqvist, in \cite{Lindqvist}, describes the first and the second eigenvalues for the $p$-laplacian. We would like to obtain other invariants which might provide us with some additional information relating the geometry of the manifold and its spectral structure.
An interesting spectral invariant on $M$ associated to the Laplace-Beltrami operator is the essential spectrum of $M$, which consists of points of the spectrum of $\triangle$ which are either accumulation points on points on the spectrum or which correspond to discrete eigenvalues of $\triangle$ with infinite multiplicity, and the greatest lower bound of the essential spectrum, $\lambda_1^{ess}(M)$. In particular, when $M$ is compact, the essential spectrum is empty and it holds the following properties $\lambda_1(M)\leq \lambda_1^{ess}(M)$ and $\lambda_1^{ess}(M)=\lim_K\lambda_1(M-K)$, where $K$ runs through all compact subsets of $M$ (see \cite{Brooks}).
Due to the difficulties in the understanding the spectrum of the $p$-Laplace operator, we shall define its {\it essential $p$-first eigenvalue}, as
\begin{equation}
\mu_{1,p}^{ess}(M):=\lim_{K} \mu_{1,p}(M-K)
\end{equation}
Where $K$ runs through all compact subsets $K$ of $M$ (More generally, we can define in a similar way its essential $p$-kth eigenvalue). With respect to essential spectrum we prove if $\theta(M)$ is the exponential volume growth of $M$ defined by
\begin{equation}
\theta(M)=\limsup_{r\to \infty} \frac1r \log(V_r(x_0))
\end{equation}
where $V_r(x_0)$ is the volume of the geodesic ball $B_r(x_0)$, then we get a Brooks-type theorem (see \cite{Brooks})

\vspace{.5cm}

\begin{thm}
If the volume of $M$ is infinity, then $\mu_{1,p}^{ess}(M)\leq \frac{\theta(M)^p}{p^p}$.
\end{thm}

\vspace{.5cm}

\section{Geometric estimates}

Following closely \cite{Bessa-Montenegro1} and \cite{Bessa-Montenegro2} we shall introduce geometric invariants associated to certain spaces of vector fields that will be used to give lower bounds for the fundamental tone for $p$-laplacian. In this direction we begin with
\begin{defi}
Let $\om \subset M$ be a domain with compact closure in a smooth Riemannian manifold $(M^n,g)$. Let $\mathfrak{X} (\om)$ be the set of all smooth vector fields, $X$, on $\om$ with sup norm $\|X\|_\infty =\sup_\om \|X\|<\infty$ (where $\|X\|=g(X,X)^{1/2}$) and $\inf_\om \div X>0 $. Define $c(\om)$ by
\begin{equation}
c(\om):= \sup\left\{ \frac{\inf_\om \div X}{\|X\|_\infty} ;\,\,\, X \in \mathfrak{X}(\om)\right\} \label{fundamental constant}
\end{equation}
\end{defi}
As remarked in \cite{Bessa-Montenegro1},  $\mathfrak{X}(\om)$ is a nonvoid set of smooth vector fields on $\om$

\noindent{\bf Proof of Theorem \ref{main1}}

Let $X\in \mathfrak{X}(\om)$ a smooth vector field and $f\in {\cal C}_0^\infty(\om)$ any positive function, then the vector field $f^pX$ has compact support in $\om$. We compute the divergence of $f^pX$

\begin{eqnarray*}
0=\int_\om \div(f^pX) \,\,d\om &=& \int_\om\Big\{ <\nabla (f^p),X> +f^p\div(X) \Big\} \,\, d\om\\
&=& \int_\om\Big\{ pf^{p-1}<\nabla f,X> +f^p\div(X) \Big\} \,\, d\om \\
&\ge & \int_\om\Big\{ -p|f|^{p-1}\|\nabla f\| \|X\| +f^p\div(X) \Big\} \,\, d\om \\
\end{eqnarray*}
by the Cauchy-Schwartz inequality. That is
\begin{equation}
0\ge \int_\om\Big\{ -p|f|^{p-1}\|\nabla f\| \|X\| +f^p\div(X) \Big\}
\,\, d\om \label{ineq1}
\end{equation}
Now, the Young inequality, for any $\alpha$, $\beta>0$
\begin{equation}
\alpha \beta \leq \frac{\alpha^p}{p}+ \frac{\beta^q}{q}, \qquad \mbox{if}\qquad  \frac1p+\frac1q=1
\end{equation}
implies that for any $\varepsilon>0$ the next inequality holds:
\begin{equation}
\alpha \beta \leq \frac{\alpha^p}{p\, \varepsilon^p}+ \frac{\varepsilon^q\beta^q}{q}\label{epsilon-young}
\end{equation}
Apply the Young inequality (\ref{epsilon-young}) to the inequality (\ref{ineq1}), letting
$$
\alpha :=p\|\nabla f\| \qquad \mbox{and} \qquad \beta :=|f|^{p-1}\|X\|
$$
to get:
\begin{eqnarray*}
0&\ge& \int_\om\bigg\{ -\frac{(p\|\nabla f\|)^p}{p \varepsilon^p} -\frac{\varepsilon^q (|f|^{p-1}\|X\|)^q}{q}  +f^p\div(X) \bigg\} \,\, d\om \\
&=& \int_\om\bigg\{ -\frac{p^{p-1}}{\varepsilon^p}\|\nabla f\|^p - \frac{\varepsilon^q\|X\|^q}q |f|^{(p-1)q}  +f^p\div(X) \bigg\} \,\, d\om \\
&=& \int_\om\bigg\{ -\frac{p^{p-1}}{\varepsilon^p}\|\nabla f\|^p + \Big(\div(X) -\frac{\varepsilon^q\|X\|^q}q\Big) |f|^p \bigg\} \,\, d\om \\
&\ge& -\frac{p^{p-1}}{\varepsilon^p}\int_\om \|\nabla f\|^p \,\, d\om +  \Big(\inf_{\om}\div(X) -\frac{\varepsilon^q}q\sup_{\om}\|X\|^q\Big) \int_\om|f|^p \,\, d\om
\end{eqnarray*}
that is
\begin{equation}
\frac{p^{p-1}}{\varepsilon^p}\int_\om \|\nabla f\|^p \,\, d\om\ge \Big(\inf_{\om}\div(X) -\frac{\varepsilon^q}q\sup_{\om}\|X\|^q\Big) \int_\om|f|^p \,\, d\om
\end{equation}
or else
\begin{eqnarray}
\int_\om \|\nabla f\|^p \,\, d\om &\ge& \frac{\varepsilon^p}{p^{p-1}} \Big(\inf_{\om}\div(X) -\frac{\varepsilon^q}q\sup_{\om}\|X\|^q\Big) \int_\om|f|^p \,\, d\om  \label{p-estimate
1}
\end{eqnarray}
Remark that when one has $\div(X)\leq 0$ on $\om$, the previous
inequality is trivial and does not bring any interesting
information. So we shall assume tacitly that $\div (X)\ge 0$ on
$\om$. Consider the function
$$
\psi(\varepsilon)=\varepsilon^p(A-B\varepsilon^q)
$$
with $A\ge 0$ and $B>0$. We will look for the maximum this function
assumes as a function of $A$ and $B$. This is a straightforward
calculation:
\begin{itemize}
\item $\psi'(\varepsilon)=\varepsilon^{p-1}[pA-(p+q)B\varepsilon^q]$
\item thus the zeroes of $\psi$ are given by
$$
\varepsilon_1 =0\qquad \mbox{and}\qquad \varepsilon_2 =
\left(\frac{pA}{(p+q)B}\right)^{1/q}
$$
\item $\psi''(\varepsilon)=\varepsilon^{p-2}[p(p-1)A-(p+q)(p+q-1)B\varepsilon^q]$
\item calculating $\psi''$ on both $\varepsilon_1$ and
$\varepsilon_2$ we get
$$
\psi''(\varepsilon_1)=0 \qquad \mbox{and} \qquad
\psi''(\varepsilon_2)=-pq\varepsilon_2^{p-2}A\leq 0
$$
\item consequently $\varepsilon_2$ is a maximum and the maximum
value of $\psi$ is given by
$$
\psi(\varepsilon_2)=\left(\frac{pA}{(p+q)B}\right)^{p/q}\frac{qA}{p+q}
=\frac{qp^{p/q}A^p}{(p+q)^pB^{p/q}}
$$
since $1+p/q=p$
\end{itemize}
We will substitute conveniently these reasonings into the integral
estimate (\ref{p-estimate 1}) letting $A=\inf_{\om}\div(X)$ and $B=\sup_{\om}\|X\|^q/q$.
Observe that:
\begin{equation*}
\max_\varepsilon \left[\varepsilon^p\Big(\inf_{\om}\div(X)
-\frac{\varepsilon^q\sup_{\om}\|X\|^q}q\Big)\right]=
\frac{q^pp^{p/q}}{(p+q)^p}\frac{(\inf_{\om}\div(X))^p}{\sup_{\om}\|X\|^p}
\end{equation*}
and consequently
\begin{eqnarray}
\frac1{p^{p-1}}\max_\varepsilon \left[\varepsilon^p\Big(\inf_{\om}\div(X)
-\frac{\varepsilon^q\sup_{\om}\|X\|^q}q\Big)\right]
=\frac1{p^p}\left(\frac{\inf_{\om}\div(X)}{\sup_{\om}\|X\|}\right)^p
\label{optimal estimate 1}
\end{eqnarray}
inserting the estimate (\ref{optimal estimate 1}) in
(\ref{p-estimate 1}) we get
\begin{eqnarray*}
\int_\om \|\nabla f\|^p \,\, d\om &\ge&
\frac1{p^p}\left(\frac{\inf_{\om}\div(X)}{\|X\|_\infty }\right)^p\int_\om |f|^p \,\, d\om
\\
&\ge& \frac1{p^p}
\left(\sup_{X\in \mathfrak{X}(\om)} \frac{\inf_{\om}\div(X)}{\|X\|_\infty}\right)^p\int_\om |f|^p \,\, d\om
\end{eqnarray*}
and thus
\begin{equation}
\int_\om \|\nabla f\|^p \,\, d\om \ge \frac{c(\om)^p}{p^p}\int_\om |f|^p \,\, d\om
\end{equation}
leading to the estimate for the fundamental tone
\begin{eqnarray}
\mu^*_{p}(\om)&=&\inf\left\{ \frac{\int_\om \|\nabla
f\|^p}{\int_\om |f|^p}:\,\,\, f\in W_+^{1,p}(\om)\setminus \{ 0\}
\right\}\ge \frac{c(\om)^p}{p^p}
\end{eqnarray}
This concludes the proof. $\square$

\vspace{.5cm}

\noindent As for a first application we prove McKean's generalized theorem \ref{mckean}:

\noindent Take for vector field $X=\nabla \rho$, the gradient of distance function from a point $o$ and observe that $\|\nabla \rho\|=1$. On the other hand $\div(\nabla \rho)=\triangle \rho$. Now, since $K_M\leq -c^2<0$ the laplacian comparison theorem implies that $\triangle \rho\ge (n-1)c$. Hence
$$\frac{(n-1)^pc^p}{p^p}\leq \frac1{p^p}\left(\frac{\div(\nabla \rho)}{\|\nabla \rho\|}\right)^p\leq \frac{c(M)^p}{p^p} \leq \mu_p^*(M)$$
concluding the proof. $\square$

\vspace{.5cm}

\noindent Now we shall obtain estimates of the fundamental tone of manifolds with curvature bounded from above:

Let $o\in M$ be a fixed point of $M$ and $\rho :B_R(o)\subset M\to\br$, the distance function from $o$, i.e., $\rho(x)={\rm dist}(o,x)$ which is smooth outside the cut locus of $o$. Set
\begin{equation}
f=\rho ^q=\rho^{p/(p-1)}
\qquad\mbox{and}\qquad
X=\|\nabla f\|^{p-2}\nabla f
\end{equation}
we determine explicitly this vector field, $X$. Since $\nabla f=q\rho^{q-1}\nabla \rho $, and $\|\nabla f\|=q\rho^{q-1}$ it follows that
\begin{eqnarray*}
X&=&(q\rho^{q-1})^{p-2}q\rho^{q-1}\nabla \rho=q^{(p-2)+1}\rho^{(q-1)(p-2)+(q-1)}\nabla \rho\\
&=&q^{p-1}\rho^{(q-1)(p-1)}\nabla \rho=q^{p-1}\rho\nabla \rho
\end{eqnarray*}

thus implying that $\|X\|\leq q^{p-1}R$ that is, $X$ is a bounded vector field on the ball $B_R(o)$. With respect to the divergency of $X$ it is immediate that
\begin{eqnarray*}
\div X&=& \div (q^{p-1}\rho\nabla \rho)=q^{p-1}\Big(<\nabla\rho,\nabla \rho> +\rho \triangle \rho\Big)\\
&=& q^{p-1}\Big(1 +\rho \triangle \rho\Big)
\end{eqnarray*}
Now, we compare the laplacian on $M$ with the laplacian on the spherically symmetric space
$$\br^n_f(\br^n, g_f=dr^2+f^2(r)d\theta^2)$$
where $g_f$ is a complete Riemannian metric given in polar coordinates - $d\theta^2$ represents the standard metric on the $(n-1)$-dimensional unit sphere in the euclidian space and $r(x)={\rm dist}(0,x)$. For $x=r\theta$, $r>0$ and $\theta\in S^{n-1}$, we remember that the hessian and laplacian of $r$ are given by
\begin{equation*}
{\rm Hess} r(x)=\frac{f'(r)}{f(r)}(g_f-dr\otimes dr)\qquad \triangle r=(n-1)\frac{f'(r)}{f(r)}
\end{equation*}
Let $K_f^{rad}(p)=-\dfrac{f''(r)}{f(r)}$ denote the radial sectional curvatures (that is the curvatures along geodesic rays from the base point $0$, containing the tangent direction to the ray) of $\br^n_f$ at $p=r\theta$. We remember the comparison theorem:

\begin{lema}[Laplacian Comparison Theorem]
Let $M$ be a complete Riemannian manifold and $x_0, x_1\in M$. Let $\gamma :[0,\rho(x_1)]\to M$ be a minimizing geodesic joining $x_0$ and $x_1$, where $\rho (x)$ is the distance function ${\rm dist}(x_0,x)$. Let $K_M^{rad}\ge $ denote the radial sectional curvatures of $M$ along $\gamma$ and assume that $K_M^{rad}(x)\leq K_f^{rad}(\rho(x))$ then
\begin{equation}
\triangle \rho(x) \ge (n-1)\frac{f'(\rho(x))}{f(\rho(x))}
\end{equation}
In particular, when $k$ is a constant
\begin{equation}
f(r)=f_k(r)=   \left\{\begin{array}{cccc}
\frac1{\sqrt{k}}\sin{\sqrt{k}r}& if &  k>0 & \\
r& if &  k=0 & \\
\frac1{\sqrt{-k}}\sinh{\sqrt{-k}r}& if &  k<0 & \\,&
\end{array}\right.
\end{equation}
turns $\br^n_{f_k}$ into the space form of constant curvature $k$. In these cases $f'(r)/f(r)$,  are given by
\begin{equation}
\frac{f'(r)}{f(r)}=\left\{\begin{array}{cccc}
\sqrt{-k} \coth (\sqrt{-k}r)& if & k<0 & \\
\frac{1}{r }& if & k=0&\\
 \sqrt{k}\cot (\sqrt{k}r)& if & k>0 ,&\quad \mbox{and $r <\pi /2\sqrt{k}$}
\end{array}\right.
\end{equation}
\end{lema}
hence, setting $\mu_f(r)=f'(r)/f(r)$:
\begin{eqnarray*}
  \div(X) \ge  q^{p-1}\Big(1 + (n-1)\rho(x)\mu_f(\rho(x))\Big)
\end{eqnarray*}
thus
\begin{equation*}
\inf_{B_R(o)}\div(X)\ge q^{p-1}\Big(1 + (n-1)\inf_{B_R(o)} \rho(x)\mu_f(\rho(x))\Big)
\end{equation*}
consequently we obtain the following estimate:
\begin{eqnarray}
\mu^*_{p}(B_R(o))\ge \left(\sup_X  \frac{\inf_{B_R(o)}\div(X)}{p\|X\|_{\chi_{\infty} (B_R(o))}|}\right)^p\ge
\left(\frac{(1 + (n-1)\inf_{B_R(o)} \rho(x)\mu_f(\rho(x)))}{pR}\right)^p
\end{eqnarray}
Summarizing, we get the following generalization of theorem (4.1) of \cite{Bessa-Montenegro1}
\begin{thm}
Let $M^n$ an $n$-dimensional complete Riemannian manifold and $B_M(q,r)$ a geodesic ball with radius $r<inj (q)$. Let $\kappa (q,r)=\sup\{ K_M(x); x\in B_M(q,r)\}$ where $K_M(x)$ are the sectional curvatures of $M$ at $x$. We have for $k>0$ that
\begin{eqnarray}
\mu^*_{p}(B_R(o))\ge \left(\sup_X  \frac{\inf_{B_R(o)}\div(X)}{p\|X\|_{\chi_{\infty} (B_R(o))}|}\right)^p\ge
\left(\frac{(1 + (n-1)\inf_{B_R(o)} \rho(x)\mu_1(\rho(x)))}{pR}\right)^p
\end{eqnarray}
\end{thm}

\noindent We now prove Theorem \ref{main2} which constitute a valuable tool for estimating the fundamental tone on open manifolds.

\vspace{.5cm}

\noindent{\bf Proof of Theorem \ref{main2}}
Let $X\in \W^{1,1}(M)$ and $f\in \cc^{\infty}_0$ In case of closed manifold or a bounded domain, $\om\subset M$ we get the estimate for any nonnegative test function $f\in \cc^{\infty}_0$ and any smooth vector field $X$ such that $supp(f^pX)\subset \subset M$:
\begin{eqnarray*}
0&=&\int_M \div(f^pX) \,\, dM=\int_M \Big[<\nabla(f^p),X> +f^p\div(X)\Big] \,\, dM\\
&=&\int_M \Big[pf^{p-1}<\nabla f,X> +f^p\div(X)\Big] \,\, dM\ge\int_M \Big[-pf^{p-1}\|\nabla f\|\|X\| +f^p\div(X)\Big] \,\, dM\\
&\ge& \int_M \left[-p\Big(\frac{\|\nabla f\|^p}p+\frac{f^{(p-1)q}\|X\|^q}q\Big) +f^p\div(X)\right] \,\, dM\\
&\ge& -\int_M \|\nabla f\|^p \,\, dM +\int_M \left(-\frac{p}q\|X\|^q+\div(X)\right)f^p \,\, dM\\
&\ge& -\int_M \|\nabla f\|^p \,\, dM + \inf_M\Big((1-p)\|X\|^q+\div(X)\Big)\int_M f^p \,\, dM
\end{eqnarray*}
where we have used the Young inequality $ab\leq \dfrac{a^p}p+\dfrac{b^q}q$ for the pair $a=\|\nabla f\|$ and $b=f^{p-1}\|X\|$ and the fact that the exponents $p,q$ are conjugate, that is $(p-1)q=p$. Thus we have
\begin{eqnarray*}
\int_M \|\nabla f\|^p \,\, dM \ge \inf_M\Big((1-p)\|X\|^q+\div(X)\Big)\int_M f^p \,\, dM
\end{eqnarray*}
or
\begin{eqnarray*}
\frac{\int_M \|\nabla f\|^p \,\, dM}{\int_M f^p \,\, dM} \ge \inf_M\Big((1-p)\|X\|^q+\div(X)\Big)
\end{eqnarray*}
for any vector field $X\in \W^{1,1}(M)$, hence we obtain
\begin{eqnarray}
\frac{\int_M \|\nabla f\|^p \,\, dM}{\int_M f^p \,\, dM} \ge \sup_{X\in \W^{1,1}(M)} \inf_M\Big((1-p)\|X\|^q+\div(X)\Big) \label{variational2}
\end{eqnarray}
and taking the greatest lower bound over all test functions $f$ on the left side  of the equation \ref{variational2} we get
\begin{eqnarray}
\mu_{p}^*(M)=\inf_{W^{1,p}_0(M)} \frac{\int_M \|\nabla f\|^p \,\, dM}{\int_M |f|^p \,\, dM} \ge \sup_{X\in \W^{1,1}(M)} \inf_M\Big((1-p)\|X\|^q+\div(X)\Big) \label{variational2}
\end{eqnarray}
This concludes the proof of the lema. $\square$

\vspace{.5cm}

Now let $u$ be the first eigenfunction associated to the eigenvalue $\mu_{1,p}$, that is
$$
\triangle_pu=\mu_{1,p}(\om)|u|^{p-2}u
$$
and consider the vector field
\begin{equation}
X=-\frac{\|\nabla u\|^{p-2}\nabla u}{|u|^{p-2}u}
\end{equation}
we calculate its norm and divergence:
\begin{equation*}
\| X\|^q=\left(\frac{\|\nabla u\|^{p-2}}{|u|^{p-1}} \|\nabla u\|\right)^q=\left(\frac{\|\nabla u \|^{p-1}}{|u|^{p-1}}\right)^q=\frac{\|\nabla u\|^{(p-1)q}}{|u|^{(p-1)q}}=\frac{\|\nabla u\|^p}{|u|^p}
\end{equation*}
and
\begin{eqnarray*}
\div(X)&=&-\div\left( \frac{\|\nabla u\|^{p-2}\nabla u}{u^{p-1}} \right)= -\frac{\div(\|\nabla u\|^{p-2}\nabla u)}{u^{p-1}} - <\|\nabla u\|^{p-2}\nabla u, \nabla \frac1{u^{(p-1)}}>\\
&=&\frac{\triangle_p  u}{|u|^{p-2}u}  +(p-1)u^{-p}\|\nabla u\|^{p-2}<\nabla u, \nabla u>
=\mu_{1,p}(\om)  +(p-1)\frac{\|\nabla u\|^p}{u^p}
\end{eqnarray*}
gathering these results
\begin{eqnarray*}
(1-p)\|X\|^q+\div(X)=(1-p)\frac{\|\nabla u\|^p}{|u|^p}+\mu_{1,p}(\om)  +(p-1)\frac{\|\nabla u\|^p}{u^p}=\mu_{1,p}(\om)
\end{eqnarray*}

We now consider estimates on the essential spectrum of $\triangle_p$. This section is based on paper of Robert Brooks \cite{Brooks}

Pick a point $x_0 \in M^n$, and for each $r>0$, denote $B_r(x_0)$ the ball centered at $p$, of radius $r$ and  $V_r(x_0)$ the volume of this ball. If one sets
\begin{equation}
\theta(M)=\limsup_{r\to \infty} \frac1r \log(V_r(x_0))
\end{equation}
it is not hard to verify, via triangle inequality, that this number does not depend on $x_0$, and is actually an invariant of the manifold, $M$ called the exponential growth of $M$. For example If $M=H^n(-c^2)$ is the hyperbolic space of constant curvature $-c^2$, then the volume of any ball is given by
$$
V(r)=\frac{\alpha(n)}{c^{n-1}}\int_0^r \sinh^{n-1}{ct}\,\, dt
$$
implying that
\begin{eqnarray*}
\theta(M) &=&\limsup_{r\to \infty} \frac1r \log\left(\frac{\alpha(n)}{c^{n-1}}\int_0^r \sinh^{n-1}(ct)\,\, dt\right)\\
&=& \lim_{r\to \infty} \frac{\sinh^{n-1}{cr}}{\int_0^r \sinh^{n-1}(ct)\,\, dt}\\
&=&(n-1)c\lim_{r\to \infty} \frac{\sinh^{n-2}{cr}\cosh{cr}}{\sinh^{n-1}(cr)}=(n-1)c
\end{eqnarray*}
in particular if $M$ is an open manifold with curvature bounded above by $\sup K_M\leq -c^2$ then, the Bishop-Gromov volume comparison gives for any $q\in M$ and any $r>0$, $vol_{M}(B_r(q))\ge vol_{H^n(-c^2)}(B_r)$, where $B_r$ is any ball of radius $r$ in $H^n(-c^2)$. Consequently, we get
  \begin{eqnarray*}
\theta (M)\ge (n-1)c=\theta (H^n(-c^2))
\end{eqnarray*}

\begin{thm}\label{essencialspec}
If the volume of $M$ is infinity, then $\mu_p^{ess}\leq \frac{\theta(M)^p}{p^p}$
\end{thm}

The proof of this theorem is, amazingly, a straightforward consequence of the lemma below (see Theorem 2 of \cite{Brooks}):

\begin{lema}\label{lemafundamental}
Let $K$ be a compact set (possibly empty) subset of $M$, and $\mu_{0,p}(M-K)=$ the greatest lower bound of the spectrum of $\triangle_p$ on $L^p(M-K)$, with Dirichlet boundary condition on $\partial$. Let $\rho(x)={\rm dist}(x_0,x)$ denote the distance function from a fixed point $x_0\in M$. If
$$
\int_{M-K}e^{-2p\alpha \rho(x)}\,\, dx \qquad \mbox{for some $\alpha$ satisfying $0<\alpha <\sqrt[p]{\mu_{0,p}(M-K)}$}
$$
then
$$
\int_{M-K}e^{2p\alpha \rho(x)}\,\, dx<\infty
$$
\end{lema}
We proceed the proof of the lema:

Consider a test function defined by $f(x)=e^{h(x)}\chi(x)$, where we assume $\chi(x)$ has compact support in $M-K$. In the sequel we will denote for simplicity $\mu_0:=\mu_{0,p}(M-K)$. Then, since
\begin{equation}
\int_{M-K}\|\nabla f\|^p \,\, dx\ge \mu_0\int_{M-K} |f|^p \,\, dx \label{aa}
\end{equation}
it follows that
\begin{eqnarray*}
\int_{M-K}\|\nabla f\|^p\,\, dx &=&\int_{M-K} e^{ph(x)}\| \chi \nabla h +\nabla \chi\|^p\,\, dx \\
&=& \int_{M-K} e^{ph(x)}|< \chi \nabla h +\nabla \chi,\chi \nabla h +\nabla \chi>|^{p/2}\,\, dx\\
&\leq&\int_{M-K} e^{ph(x)}| \chi^2 \|\nabla h\|^2 +2\chi |<\nabla h,\nabla \chi>|+\|\nabla \chi\|^2 |^{p/2}\,\, dx\\
&\leq&\int_{M-K} |f|^p\|\nabla h\|^p+ e^{ph(x)}\Big[ 2^{p/2}\chi^{p/2} \|\nabla h\|^{p/2}\|\nabla \chi\|^{p/2}+\|\nabla \chi\|^p \Big]\,\, dx
\end{eqnarray*}
that is
\begin{eqnarray}
\int_{M-K}\|\nabla f\|^p\,\, dx &\leq&\int_{M-K} |f|^p\|\nabla h\|^p\,\, dx \nonumber\\
&&+ \int_{M-K}e^{ph(x)}\Big[ 2^{p/2}\chi^{p/2} \|\nabla h\|^{p/2}\|\nabla \chi\|^{p/2}+\|\nabla \chi\|^p \Big]\,\, dx \label{bb}
\end{eqnarray}
hence, comparing (\ref{aa}) with (\ref{bb}) we get
\begin{eqnarray}
 \int_{M-K} |f|^p\Big( \mu_0-\|\nabla h\|^p\Big) \,\, dx \leq  \int_{M-K}e^{ph(x)}\Big[ 2^{p/2}\chi^{p/2} \|\nabla h\|^{p/2}\|\nabla \chi\|^{p/2}+\|\nabla \chi\|^p \Big]\,\, dx \label{cc}
\end{eqnarray}

We now suppose that $\|\nabla h\|\leq \alpha <\mu_o^{1/p}$, and, for some exhaustion $K_1\subset \ldots\subset K_n\subset \ldots \subset M-K$ through compact subsets of $M-K$ we set, for each $i\in \bn$ and any fixed positive number $d$
\begin{equation}
\chi_{i,d} (x)=\left\{ \begin{array}{ccc}
0,&\mbox{if}& x\in M-(K\cup K_i)\\
\dfrac{1}{d}\rho(x,M-K_i),&\mbox{if}& 0\leq \rho(x,M-K_i)\leq d\\
1,&\mbox{if}& \rho(x,M-K_i)\ge d
\end{array}\right.
\end{equation}
Then it is immediate to verify that $\|\nabla \chi_{i,d}\|\leq \dfrac1d$, and $\nabla \chi_{i,d}$ is supported in the tubular neighborhood $T_d(\partial K_i)$ of radius $d$ around $\partial K_i$. Observe that according $i\to \infty$ one gets $\partial K_i\to \partial (M-K)=\partial K$ thus $T_d(\partial K_i)\to T_d(\partial K)$.

We now suppose that $\|\nabla h\|\leq \alpha <\mu_o^{1/p}$ then, taking for each $i$ the test function $f_i(x)=e^{h(x)}\chi_{i,d}(x)$ from (\ref{cc}) and the above restriction on $\|\nabla h\|$ we obtain
\begin{eqnarray*}
 \int_{M-K} |f_i|^p\Big( \mu_0-\alpha^p\Big) \,\, dx &\leq& \int_{M-K} |f_i|^p\Big( \mu_0-\|\nabla h\|^p\Big) \,\, dx \\
&\leq& \int_{T_d(\partial K_i)}e^{ph(x)}\Big[ 2^{p/2}\chi_i^{p/2} \|\nabla h\|^{p/2}\|\nabla \chi_i\|^{p/2}+\|\nabla \chi_i\|^p \Big]\,\, dx\\
&\leq& \int_{T_d(\partial K_i)}e^{ph(x)}\Big[ 2^{p/2} \|\nabla h\|^{p/2}\frac1{d^{p/2}}+\frac1{d^p} \Big]\,\, dx\\
&\leq&  \left[ \mu_o^{1/2}\left(\frac2{d}\right)^{p/2} +\frac1{d^p} \right]\int_{T_d(\partial K_i)}e^{ph(x)}\,\, dx
 \end{eqnarray*}
Consequently, for $e^{ph}\in L^1(M-K)$, letting $i\to \infty$ we get
\begin{eqnarray*}
 \int_{M-K} |f|^p\Big( \mu_0-\alpha^p\Big) \,\, dx \leq \left[ \mu_o^{1/2}\left(\frac2{d}\right)^{p/2} +\frac1{d^p} \right]\int_{T_d(\partial K)}e^{ph(x)}\,\, dx
 \end{eqnarray*}

Under the assumption $\int_{M-K} e^{-p\alpha \rho(x)}\,\, dx<\infty$ we may define a sequence of functions $h_j$ as follows
$$h_j(x)=\min\{ \alpha \rho(x),j-\alpha\rho(x)\}$$
 Observe that for each $j$, $\|\nabla h_j\|\leq \alpha$, and that this is an increasing sequence converging pointwisely to  $h(x)=\alpha \rho(x)$. Thus, for $j$ sufficiently large (for instance $j>2\alpha.{\rm dist}(x_0,T_d(\partial K))$), we have
\begin{eqnarray*}
 \int_{M-K} e^{ph_j(x)}\Big( \mu_0-\alpha^p\Big) \,\, dx \leq \left[ \mu_o^{1/2}\left(\frac2{d}\right)^{p/2} +\frac1{d^p} \right]\int_{T_d(\partial K)}e^{p\alpha\rho(x)}\,\, dx
 \end{eqnarray*}
 so
 $$
 \int_{M-K} e^{ph_j(x)}\,\, dx \leq C
 $$
where $C$ is a finite constant, independent of $j$, given by
$$C=C(\mu_0,d,\alpha,p,|T_d(\partial K)|)= \left(\frac{\mu_o^{1/2}(2d)^{p/2} +1}{d^p(\mu_0-\alpha^p)}\right) |T_d(\partial K)|e^{p\alpha.{\rm dist}(x_0,T_d(\partial K))}$$
and $|T_d(\partial K)|=volume(T_d(\partial K))$.
Taking the limit as $j\to \infty$ gives
\begin{equation}
\int_{M-K}e^{p\alpha\rho(x)}\,\, dx\leq C
\end{equation}
This concludes the theorem. $\square$

Now we can proceed the

\noindent{\bf Proof of Theorem \ref{essencialspec}}

Since $M$ has exponential volume growth it follows that $r\mapsto \dfrac{\log(V(r))}{r}$ converges, for $r\to \infty$, to $\theta$. Thus if $p\alpha >\theta $, and $\delta>0$ is any positive constant such that $p\alpha >\theta+\delta$ take $\epsilon =p\alpha -\theta-\delta$, then this implies that there exists $M >0$ such that for $r>M$ it holds $V(r)<e^{(\theta+\epsilon)r}=e^{(p\alpha-\delta)r}$

\begin{eqnarray*}
\int_{M-K}e^{-p\alpha\rho(x)}\,\, dx&\leq& \sum_{r=1}^\infty [V(r)-V(r-1)]e^{-p\alpha (r-1)}\\
&\leq&\sum_{r=1}^\infty V(r)e^{-p\alpha (r-1)}[e^{p\alpha} -1]\\
&<&\sum_{r=1}^M V(r)e^{-p\alpha (r-1)}[e^{p\alpha} -1]+\sum_{r=M}^\infty e^{-\delta r}[e^{p\alpha }(e^{p\alpha} -1)]<\infty
\end{eqnarray*}
Since the last expression is a sum of a geometric series.

We now conclude from lemma \ref{lemafundamental} that if $p\alpha >\theta$ and $\alpha <\sqrt[p]{\mu_{p,1}(M-K)}$ we have $\int_{M-K}e^{p\alpha\rho(x)}\,\, dx<\infty$. But this is impossible since $M$ has infinite volume.

Therefore, there is no such $\alpha$, and we have shown $\mu_{1,p}(M-K)\leq \dfrac{\theta(M)^p}{p^p}$. Taking the limit over arbitrarily large $K$, we see that $\mu_{p}^{ess}(M)\leq \dfrac{\theta(M)^p}{p^p}$.  $\square$

Now we point some consequences (see \cite{Brooks}).

\begin{cor}
If $M$ has subexponential growth then $\mu_{1,p}^{ess}(M)=0$
\end{cor}
In particular, when $M^n$ is a manifold of polynomial growth one has $\mu_{1,p}^{ess}(M)=0$

\noindent{\bf Proof} If $M$ has subexponential growth then $\theta(M) =0$, thus $0\leq \mu_{p}^{ess}(M)\leq 0$. $\square$

To state next corollary we need Cheeger isoperimetric constant, $h(M)$:
\begin{equation}
h(M)=\inf_N \frac{area(N)}{vol({\rm int} N)}
\end{equation}
where $N$ runs over all compact $(n-1)$-dimensional submanifolds of $M$ dividing $M$ into two components, and ${\rm int}N$ denotes the bounded component.

\begin{cor}
$h(M)\leq \theta(M)$. Furthermore, equality holds if the ratio $S(r)/V(r)$ tends to $h(M)$ as $r\to\infty$, where $S(r)$ denotes the surface area of the distance sphere of radius $r$. In this case $\mu_{1,p}(M)=\dfrac{\theta(M)^p}{p^p}$
\end{cor}

\end{document}